\magnification=\magstep1
\hsize16truecm
\vsize23.5truecm
\topskip=1truecm
\raggedbottom
\abovedisplayskip=3mm
\belowdisplayskip=3mm
\abovedisplayshortskip=0mm
\belowdisplayshortskip=2mm
\normalbaselineskip=12pt
\normalbaselines
\input amssym.def
\input amssym.tex
\font\titlefont= cmcsc10 at 12pt
\def\F{\Bbb F}

\def\Z{\Bbb Z}

\def\P{\Bbb P}

 

%
\catcode`\@=11
\font\tenmsa=msam10
\font\sevenmsa=msam7
\font\fivemsa=msam5
\font\tenmsb=msbm10
\font\sevenmsb=msbm7
\font\fivemsb=msbm5
\newfam\msafam
\newfam\msbfam
\textfont\msafam=\tenmsa  \scriptfont\msafam=\sevenmsa
  \scriptscriptfont\msafam=\fivemsa
\textfont\msbfam=\tenmsb  \scriptfont\msbfam=\sevenmsb
  \scriptscriptfont\msbfam=\fivemsb
\def\hexnumber@#1{\ifcase#1 0\or1\or2\or3\or4\or5\or6\or7\or8\or9\or
	A\or B\or C\or D\or E\or F\fi }
\edef\msa@{\hexnumber@\msafam}
\edef\msb@{\hexnumber@\msbfam}
\mathchardef\square="0\msa@03
\mathchardef\subsetneq="3\msb@28
\mathchardef\ltimes="2\msb@6E
\mathchardef\rtimes="2\msb@6F
\def\Bbb{\ifmmode\let\next\Bbb@\else
	\def\next{\errmessage{Use \string\Bbb\space only in math mode}}\fi\next}
\def\Bbb@#1{{\Bbb@@{#1}}}
\def\Bbb@@#1{\fam\msbfam#1}
\catcode`\@=12
%
%
%
\def\vandaag{\number\day\space\ifcase\month\or
 januari\or februari\or  maaart\or  april\or mei\or juni\or  juli\or
 augustus\or  september\or  oktober\or november\or  december\or\fi,
\number\year}
\def\today{\ifcase\month\or
 Jan\or Febr\or  Mar\or  Apr\or May\or Jun\or  Jul\or
 Aug\or  Sep\or  Oct\or Nov\or  Dec\or\fi
 \space\number\day, \number\year}
\vskip 6.5pc
\noindent
\font\eighteenbf=cmbx10 scaled\magstep2
\font\titlefont=cmcsc10 at 10pt
\vskip 2.0pc
\centerline{\eighteenbf  Constructing Curves over Finite Fields}
\bigskip
\centerline{\eighteenbf with Many Points }
\bigskip
\centerline{\eighteenbf by Solving Linear Equations}
\vskip 2pc
\centerline{\titlefont Gerard van der Geer \& Marcel van der Vlugt}
\bigskip
The purpose of this note is to exhibit an elementary method to construct
explicitly curves over finite fields with many  points. Despite its
elementary character the method is very efficient and can be regarded as a
partial substitute for the use of class field theory. Many of the results
on the existence of curves with a large number of points obtained from
class field theory or Drinfeld modules can thus be reproduced with {\sl
explicit} curves and many new examples can easily be obtained. Curves with
many points find applications in coding theory and the theory of
low-discrepancy sequences and here explicitness is often essential.
\par
\bigskip
\centerline{\bf The Method}
\bigskip
\noindent
Let $k=\F_q$ be a finite field of cardinality $q=p^m$ with $p$ prime. We
consider fibre products $C_F$ over a base curve $C$  of 
Artin-Schreier extensions $C_f \to C$  defined by function field extensions
$k(C_f)=k(C)(z)$, where $z$ satisfies an equation $z^p-z=f$ and $f$ 
runs through a basis of functions $f \in F$ with $F$  a suitable
$\F_p$-vectorspace of $k(C)$.  In order to get a curve
$C_F$ with many points we shall start with a curve $C$ over $k$ with many
rational points. We select a (large) set of distinct  rational points
${\cal P}$ on $C$ and an effective divisor $D$ defined over $k$ whose
support is disjoint from ${\cal P}$. We consider a $\F_p$-vector space  
$$
F \subseteq L(D)=\{ g \in
k(C): (g) + D\geq 0 \} \cup \{ 0 \}
$$ 
of functions $f$ satisfying
$$
\eqalign{
{\rm i)}& \quad F \cap \{ g^p-g: g \in k(C) \, \} = \{ 0 \},
\cr
{\rm ii)}& \quad{\rm Tr}_{q/p}(f(P))=0 \quad {\hbox{\rm for all } } P
\in {\cal P},\cr
}
$$
where ${\rm Tr}_{q/p}$ is the trace map from $\F_q$ to $\F_p$. The fibre
product
$C_F$ is a covering of
$C$ in which all the rational points $P \in {\cal P}$ are completely split
by ii), so that we may expect many points. The genus of $C_F$ is determined
by the genera of the building blocks  $C_f$ with $f\in F-\{ 0 \}$ and these
follow from the Hurwitz-Zeuthen formula. A sensible choice of $D$ must be
made in order to keep the genus of $C_F$ relatively small. Linear algebra 
gives estimates for the dimensions of  subspaces $F$ satisfying i) and
ii).
\eject
\bigskip
\centerline{\bf Covers of Curves}
\bigskip
\noindent
Let $C$ be a smooth complete irreducible curve defined over $k$ and suppose
that $f_1,\ldots,f_r$ is a basis of a space of functions $F$ satisfying
condition i). Each curve $C_f$ defined by the equation $z_i^p-z_i=f_i$  is
an irreducible cover of $C$; up to $C$-isomorphism the fibre product 
$$
C_F=C_{f_1} \times_C\times \ldots \times_C C_{f_r}
$$
does not depend on the chosen basis. The Jacobian
$J(C_F)$ of $C_F$ decomposes up to isogeny as product
$$
J(C)\times \prod_{f \in \P(F)}P_f,
$$
where $P_f$ is the Prym variety of $C_f \to C$ (i.e.\ the connected
component of the kernel of the norm map ${\rm Nm}: J(C_f) \to J(C)$).
The sum is over a complete set of representatives $f \in F$ modulo 
the equivalence relation $f \sim c\, f$ with $c \in \F_p^*$ (which we
write as $f \in \P(F)$).
\par
If we denote the trace of Frobenius of $C$ (resp.\ $C_f$, $C_F$) by 
$\tau$ (resp.\ $\tau_f$, $\tau_F$) then we find the relation
$$
\tau_F= \tau + \sum_{f \in \P(F)}(\tau_f-\tau)
$$
from which we can compute the number of points of $C_F$. The genus of $C_F$
is given similarly by
$$
g(C_F)= g(C) + \sum_{f \in \P(F)} (g(C_f)-g(C)).
$$
\par
\smallskip
\noindent
\centerline{\bf Curves arising from Spaces of Solutions}
\smallskip
\noindent
To guarantee that we have sufficiently many functions we use Riemann-Roch.
Let ${\cal P}=\{P_1,\ldots, P_n\}$ be the set of $k=\F_q$-rational points
on $C$. For an effective  divisor $D= \sum n_Q Q$ defined over $\F_q$ we
denote by $[D/p]$ the divisor $\sum[n_Q/p] Q$ with $[x]$ the
greatest integer function. We set
$$
\delta = \# ({\rm Supp}(D)\cap {\cal P}).
$$
\smallskip
\noindent
\proclaim Proposition 1. For a curve $C$ and a divisor $D$ as above there
exists for every integer $r$ with
$$
1 \leq r \leq \dim_{\F_p}(L(D))-\dim_{\F_p}(L([D/p])+1-n+\delta
$$
a $\F_p$-subspace $F_r \subset L(D)$ of dimension $r$ such that for $f\in
F_r-\{ 0 \}$ the curve $C_f$ satisfies
$$
\# C_f(\F_q)=p(n-\delta)+\epsilon_f,
$$
where $\epsilon_f\geq 0$ is the number of rational points of $C_f$ lying
over points in ${\rm Supp}(D)$.
\par
\noindent
{\sl Proof.} Let ${\bar{\F}}$ denote an algebraic
closure of $\F_q$. From [V] it follows that the $\F_p$-vector space 
$$
V= \{ f \in L(D) \colon f= g^p-g {\hbox { \rm for some } } g \in
{\bar{\F}}(C)
\}
$$
has dimension $\dim_{\F_p}(L([D/p]))$. Choose now a linear $\F_p$-space
$W\subset L(D) $  with $W \cap V=\{ 0 \}$ and with dimension
$\dim_{\F_p}(L(D))- \dim_{\F_p}(L([D/p]))$. Let $c \in \F_q$ with ${\rm
Tr}_{q/p}(c)\neq 0$, set
$$
\tilde W= W \oplus \F_p\cdot c
$$
and solve in
$\tilde W$  the linear equations
$$
{\rm Tr}_{q/p}(f(P_i))=0 \quad {\hbox {\rm for all }} P_i\not\in {\rm
Supp}(D).\eqno(1)
$$
The solution space of this homogeneous system of linear equations over
$\F_p$ has dimension $\geq \dim \tilde W - (n-\delta)$. The equations (1) 
express the condition that all the points $P_i$ which do not belong to
${\rm Supp}(D)$ are completely split. This gives the expression for
$\#C_f(\F_q)$. 
$\square$
\bigskip
We now form the fibre product $C_{F}$ associated to a space of
solutions $F$ of (1) and we obtain:
\smallskip
\noindent
\proclaim Corollary. With the notations as in Proposition 1 there exists
for every  $r$ with 
$$1\leq r \leq \dim_{\F_p}(L(D))-\dim_{\F_p}(L([D/p]))+1-n+\delta
$$
a $r$-dimensional subspace $F_r\subset L(D)$ such that
$$
\# C_{F_r}(\F_q)= p^r(n - \delta)+ \delta + \sum_{f \in
\P(F_r)}(\epsilon_f-\delta).
$$
\par
\noindent
{\bf Remark.} For $f \in F_r, f\neq 0$ the Hurwitz-Zeuthen formula implies:
$$
2g(C_f)-2= p(2g(C)-2)+(p-1)\sum_{Q}((-v_Q^*(f)+1)\deg(Q)),
$$
where the summation is over those $Q$ in ${\rm Supp}(D)$ with $-v_Q^*(f) >
0$. Here $v_Q^*(f)$ denotes the Artin-Schreier reduced order at the pole
of $f$ in $Q$.
\bigskip
\noindent
\centerline{\bf Examples}
\bigskip
\noindent
Examples can now be obtained in a straightforward way. By considering the
special cases
$$
\eqalign{
D=2Q\qquad & {\hbox{ \rm with $Q$ a point of degree $d$ on $C$} },\cr
D=kP\qquad & {\hbox{ \rm with $k \in \Z_{>0}$ and $P$  a rational  point of
$C$} },\cr 
D=2P+2Q\qquad & {\hbox{ \rm with $P$ and $Q$ as above}},\cr
}$$
one can reproduce, but in {\sl explicit form}, many of the examples given
in the papers of Niederreiter and Xing (see References). There the curves
were implicitly given by subfields of narrow ray class fields determined by
Drinfeld modules of rank $1$ over a curve. However, if  the base curve has
genus $\geq 1$ it is rather difficult to find explicit forms of the
curves.  The same holds for examples constructed by class field
theory in [L], [Sch] and [S1,2].
\smallskip
\noindent
{\bf Example 1.} In [S3] Serre constructed curves $C$ over $\F_2$ with
invariants  $(g(C), \#C(\F_2))$ $= (10,13), (11,14)$ and $(13,15)$ as
covers of the elliptic curve $E$ defined by $y^2+y=x^3+x$ using the group
structure of $E(\F_2)$. We show how these examples can be recovered with
our method. The rational points of $E$ are
$P_{\infty}$, $P_1=(0,0)$, $P_2=(1,0)$, $P_3=(1,1)$ and $P_4=(0,1)$.
Consider the space $\tilde W$ for  $D= 3P_{\infty} + P_1+P_3$ with
$\F_2$-basis
$$
1, x, y , x/(x+y), xy/(x+y).
$$
We now require in $\tilde W$ that $P_2$ and $P_4$ are completely split. The
solution space of (1) has basis $f_1=x + x/(x+y)$, $f_2=1+x+y$ and
$f_3=xy/(x+y)$. The $2$-dimensional $\F_2$-subspace $F=\langle f_1,
f_2+f_3\rangle$ yields a curve of genus $10$ with $13$ points. 
\par
With $D=3P_{\infty} + 3P_4$ we can take $\tilde W= \langle 1, y/x, x, y,
y/(x(x+y+1)\rangle$. Then $F= \langle 1+y/x+x+y, (xy+y^2)/x(x+y+1)\rangle$
yields the example with $(g,\#)=(11,14)$.
\par
With $D=7P_{\infty} +P_4$ and splitting $\{ P_1,P_2,P_3\}$ we find a
$3$-dimensional  space $H=\langle 1+x+y+y/x,
y+xy, y+x^2y\rangle$. The subspace $F=\langle 1+x+y+y/x, xy+x^2y\rangle$
produces a curve $C_F$ with $(g(C_F),\#C_F(\F_2))=(13,15)$. Other
$2$-dimensional spaces produce curves with
$(g,\#)=(11,13),(12,13),(13,13),(14,13),(14,14)$ and $(14,15)$. The curve
$C_H$ has $(g(C_H),\#C_H(\F_2))=(29,25)$.
\smallskip
\noindent
{\bf Example 2 i).} Consider the elliptic curve $C$ defined by
$y^3-y=x^2-1$ over $\F_3$ with $\#C(\F_3)=7$. Letting $P_ {\infty}=
(0:1:0)$, $Q=(1,0)$ and  $D=8 P_{\infty}+2Q$ we can take $\tilde W$ with
$\F_3$-basis 
$$
1,(y+2)/(2x+y+1), (x+y+1)/y, y, x+2, y^2, (x+2)y, (x+2)y^2, (x+2)^2y.
$$
We solve in $\tilde W$ for the remaining  $5$ rational points
$P_i$ ($i=1,\ldots,5$) of $C$ the equations 
$$
{\rm Tr}_{3/3}(f(P_i))= f(P_i)=0\qquad i=1,\ldots, 5
$$
and we get a $4$-dimensional solution space $F$ generated by
$$
\eqalign{ f_1={y+2 \over 2x+y+1} + y + x+y^2,& \quad 
f_2={y+2 \over 2x+y+1}+{x+y+1\over y} + x +2xy,\cr
 f_3= 2+ {y+2 \over
2x+y+1} + y+2y^2+&xy^2, \quad {\rm and} \quad f_4=2y+x^2y \cr}
$$
with the following curves
\medskip
\vbox{
\noindent{\bf }
\medskip\centerline{\def\quad{\hskip 0.6em\relax}
\def\quod{\hskip 0.5em\relax }
\vbox{\offinterlineskip
\hrule
\halign{&\vrule#&\strut\quod\hfil#\quad\cr
height2pt&\omit&&\omit&&\omit&\cr
&$f_i=$&&$g(C_{f_i})$&&$\#C_{f_i}(\F_3)$&\cr
height2pt&\omit&&\omit&&\omit&\cr
\noalign{\hrule}
height2pt&\omit&&\omit&&\omit&\cr
&$f_1$&&$9$&&$17$&\cr
&$f_2$&&$10$&&$17$&\cr
&$f_3$&&$12$&&$17$&\cr
&$f_4$&&$10$&&$19$&\cr
height2pt&\omit&&\omit&&\omit&\cr
\cr } \hrule}
}}
\medskip
\noindent
Analyzing the subspaces of $F$ we obtain the following table, in which
under $N_q(g)$ we list the best curve so far and the Oesterl\'e upper bound
(cf.\ [G-V]).
\medskip
\vbox{
\noindent{\bf }
\medskip\centerline{\def\quad{\hskip 0.6em\relax}
\def\quod{\hskip 0.5em\relax }
\vbox{\offinterlineskip
\hrule
\halign{&\vrule#&\strut\quod\hfil#\quad\cr
height2pt&\omit&&\omit&&\omit&&\omit&\cr
&$F$&&$g(C_F)$&&$\# C_F(\F_3))$&&$N_q(g)$&\cr
height2pt&\omit&&\omit&&\omit&&\omit&\cr
\noalign{\hrule}
height2pt&\omit&&\omit&&\omit&&\omit&\cr
&$\langle f_1,f_2\rangle$&&$35$&&$47$&&$[38-51]$&\cr
&$\langle f_1+2f_3,f_4\rangle$&&$36$&&$46$&&$[36-52]$&\cr 
&$\langle f_1,f_3\rangle$&&$39$&&$46$&&$[42-56]$&\cr 
&$\langle f_1,f_2,f_3\rangle$&&$128$&&$136$&&$\leq 149$&\cr
height2pt&\omit&&\omit&&\omit&&\omit&\cr
\cr } \hrule}
}}
\medskip
\noindent
{\bf Example 2 ii).} Consider again the elliptic curve $C$ given by
$y^3-y=x^2-1$ over $\F_3$ with $\#C(\F_3)=7$.  Let $D$ be the divisor $3Q$
with $Q$ a point of degree $4$ on $C$ given by the ideal
$(x^2+xy+1,x^4+x^3+2)$. This is the Frobenius orbit of $(\alpha,
\alpha^3+\alpha^9)$.  Requiring that the
$7$ rational points are completely split leads to a $2$-dimensional space
$F$ with basis
$$
1+f_2+2f_3+2f_1^2+f_2^2+f_3^2 \qquad {\rm and}\qquad
2+f_1^2+f_2^2+f_1f_2f_3,
$$
where 
$f_1=(2x+y)/(x^2+xy+1), f_2=(2x^2+1)/(x^2+xy+1),
f_3=(2x^2+y^2)/(x^2+xy+1).
$
Then $C_F$ is a curve of genus $49$ and $\#C_F(\F_3)=63$. (Best result
so far $N_q(g) \in [49-67]$.)
\smallskip
\noindent
{\bf Example 3.} Let $\F_4=\F_2(\alpha)$ and let $C$ be the elliptic curve 
defined by $y^2+y=x^3$ with $\#C(\F_4)=9$. We use $P_{\infty}=(0:1:0)$ and
set $D=11P_{\infty}$. Elementary calculations give $\dim_{\F_2}(\tilde
W)=13$ and we may take as a $\F_2$-basis:
$$\tilde W= \langle  \alpha, x,\alpha x, y, \alpha y, xy, \alpha xy, x^2y,
\alpha x^2y, y^3, \alpha y^3, xy^3, \alpha xy^3 \rangle.
$$
Splitting the $8$ rational
points different from $P_{\infty}$ yields $8$ equations in $12$
variables with  a $5$-dimensional solution space. We list properties of
the curves $C_{f_i}$ corresponding to a basis
$f_i$ ($i=1,\ldots, 5$) of the solution space  in the following table:
\medskip
\vbox{
\noindent{\bf }
\medskip\centerline{\def\quad{\hskip 0.6em\relax}
\def\quod{\hskip 0.5em\relax }
\vbox{\offinterlineskip
\hrule
\halign{&\vrule#&\strut\quod\hfil#\quad\cr
height2pt&\omit&&\omit&&\omit&\cr
&$f_i=$&&$g(C_{f_i})$&&$\#C_{f_i}(\F_4)$&\cr
height2pt&\omit&&\omit&&\omit&\cr
\noalign{\hrule}
height2pt&\omit&&\omit&&\omit&\cr
&$x+xy+x^2y$&&$5$&&$17$&\cr
&$\alpha x+\alpha xy+\alpha^2 x^2y$&&$5$&&$17$&\cr
&$y^3$&&$6$&&$17$&\cr
&$x+xy^3$&&$7$&&$17$&\cr
&$\alpha x+\alpha xy^3$&&$7$&&$17$&\cr
height2pt&\omit&&\omit&&\omit&\cr
\cr } \hrule}
}}
\medskip
\noindent
The corresponding fibre products yield the following harvest:
\medskip
\vbox{
\noindent{\bf }
\medskip\centerline{\def\quad{\hskip 0.6em\relax}
\def\quod{\hskip 0.5em\relax }
\vbox{\offinterlineskip
\hrule
\halign{&\vrule#&\strut\quod\hfil#\quad\cr
height2pt&\omit&&\omit&&\omit&&\omit&&\omit&\cr
&$F$&&$g(C_F)$&&$\# C_F(\F_4))$&&$N_q(g)$&&
{\hbox {\rm comments} }&\cr
height2pt&\omit&&\omit&&\omit&&\omit&&\omit&\cr
\noalign{\hrule}
height2pt&\omit&&\omit&&\omit&&\omit&&\omit&\cr
&$\langle f_1,f_2\rangle$&&$13$&&$33$&&$33$&&best possible&\cr
&$\langle f_1,f_3\rangle$&&$15$&&$33$&&$\leq 37$&&cf.\ Ex. 4.14 in
[N-X1]&\cr 
&$\langle f_1,f_2,f_3\rangle$&&$33$&&$65$&&$\leq 66$&&cf.\
[L]&\cr 
&$\langle f_1,f_2,f_4\rangle$&&$37$&&$65$&&$\leq 72$&&cf.\ [N-X2]&\cr
&$\langle f_1,f_3,f_4\rangle$&&$39$&&$65$&&$\leq 75$&&new&\cr
height2pt&\omit&&\omit&&\omit&&\omit&&\omit&\cr
\cr } \hrule}
}}
\medskip
\noindent
\par
By employing this method systematically one can expect a
substantial number of examples improving the tables (cf.\ [G-V]). We
refer to work in progress by V.\ Shabat.
\eject
\bigskip
\centerline{\bf References}
\smallskip 
\noindent
\smallskip \noindent
[G-V]  G.\ van der Geer, M.\ van der Vlugt : Tables for the function
$N_q(g)$. Regularly updated tables at: {\hbox { { \tt
http://www.wins.uva.nl/ $ \tilde{ }$geer }} }.
\smallskip \noindent
[L] K.\ Lauter: Ray class field constructions of curves over finite fields
with many rational points. In: {\sl Algorithmic Number Theory} (Talence
1996), H.\ Cohen Ed., Lecture Notes in Computer Science 1122, Springer,
Berlin, 1996, p.\ 187-195.
\smallskip \noindent
[N-X1] H.\ Niederreiter, C.\ P.\ Xing: Cyclotomic function fields, Hilbert
class fields and global function fields with many rational places. {\sl
Acta Arithm.\ \bf 79} (1997), p.\ 59-76.
\smallskip \noindent
[N-X2] H.\ Niederreiter, C.\ P.\ Xing: Drinfeld modules of rank 1 and
algebraic curves  with many rational points II. {\sl Acta Arithm.\ \bf 81}
(1997), p.\ 81-100.
\smallskip \noindent
[Sch] R.\ Schoof: Algebraic curves and coding theory. UTM 336, Univ.\
of Trento, 1990.
\smallskip \noindent
\smallskip \noindent
[S1] J-P.\ Serre : Sur le nombre de points rationnels d'une  courbe
alg\'ebrique sur un corps fini.  {\sl C.R.\ Acad.\ Sci.\ Paris \bf 296},
S\'erie I (1983), p. 397-402. (= Oeuvres III, No.\ 128,  p.\ 658-663).
\smallskip \noindent
[S2] J-P.\ Serre: Rational points on curves over finite fields.
Notes of lectures at Harvard University 1985.
\smallskip \noindent
[S3] J-P.\ Serre: Letter to G.\ van der Geer, 1 September 1997.
\smallskip
\noindent
[V] M.\ van der Vlugt: A new upper bound for the dimension of
trace codes. {\sl Bull.\ London Math.\ Soc.\ \bf 23} (1991), p.\
395-400.
\smallskip \noindent
[X-N]  C.\ P.\ Xing, H.\ Niederreiter: Modules de Drinfeld et courbes
alg\'ebriques ayant beaucoup de points rationnels.  {\sl C.R.\ Acad.\
Sci.\ Paris \bf 322}, S\'erie I (1996), p. 651-654.
\bigskip
\settabs3 \columns
\+G. van der Geer  &&M. van der Vlugt\cr
\+Faculteit
Wiskunde en Informatica &&Mathematisch Instituut\cr
\+Universiteit van
Amsterdam &&Rijksuniversiteit te Leiden \cr
\+Plantage Muidergracht 24&&Niels Bohrweg 1 \cr
\+1018 TV Amsterdam
&&2333 CA Leiden \cr
\+The Netherlands &&The Netherlands \cr
\+{\tt geer@wins.uva.nl} &&{\tt vlugt@wi.leidenuniv.nl} \cr

\bye